\documentclass{amsart}

\usepackage{tikz}
\usetikzlibrary{matrix,arrows}
\usepackage{geometry,latexsym,amssymb}

\theoremstyle{plain}
\newtheorem{Thm}{Theorem}

\newtheorem{Lem}{Lemma}

\newenvironment{Prf}{{\bf Proof:} }{\hfill $\Box$
\mbox{}}

\theoremstyle{definition}
\newtheorem{Def}{Definition}

\theoremstyle{remark}

\numberwithin{equation}{section}

\begin{document}
\title[On $G$-sequential connectedness]%
   {On $G$-sequential connectedness}
\author{H\"usey\.{I}n \c{C}akall\i*, and Osman Mucuk**\\ *Department of Mathematics, Maltepe University, Istanbul, Turkey \\
          ** Department of Mathematics, Erciyes University,
           Kayseri, Turkey}

\address{H\"usey\.{I}n \c{C}akall\i\\
          Maltepe University, Department of Mathematics, Marmara E\u{g}\.{I}t\.{I}m K\"oy\"u, TR 34857, Maltepe, \.{I}stanbul-Turkey \; \; \; \; \; Phone:(+90216)6261050 ext:2248, \;  fax:(+90216)6261113}

\email{hcakalli@maltepe.edu.tr; hcakalli@@gmail.com}

\address{Osman Mucuk \\
           Department of Mathematics, Erciyes University,
           Faculty of Science, Kayseri, Turkey}

\email{mucuk@erciyes.edu.tr; mucukosman@gmail.com}

\keywords{Sequences, series, summability, sequential closure, $G$-sequential continuity, $G$-sequential connectedness}
\subjclass[2000]{Primary: 40J05 ; Secondary: 54A05, 22A05}

\date{\today}

\begin{abstract}
Recently, Cakalli has introduced a concept of $G$-sequential connectedness in the sense that a non-empty subset $A$ of a  Hausdorff topological group $X$ is $G$-sequentially connected if there are no  non-empty, disjoint $G$-sequentially closed subsets $U$ and $V$  meeting $A$ such that $A\subseteq U\bigcup V$. In this paper we investigate further properties of  $G$-sequential connectedness and prove some interesting theorems.

\end{abstract}

\maketitle

\section{Introduction}
The concept of connectedness and any concept related to connectedness play a very important role not only in pure mathematics but also in other branches of science involving mathematics especially in geographic information systems, population modeling and motion planning in robotics.

In \cite{ConnorGrosse}, Connor and Grosse-Erdmann have investigated the impact of changing the definition of the convergence of sequences on the structure of sequential continuity of real functions.  Cakalli \cite{CakalliSequentialdefinitionsofcompactness} extended this concept to topological group setting and has introduced  the concepts of $G$-sequential compactness and $G$-sequential continuity; and has investigated some results in this generalized setting (see also \cite{CakalliOnGcontinuity}).

One is often relieved to find that the standard closed set definition of connectedness for topological  spaces can be replaced by a sequential definition of connectedness. That many of the properties of connectedness of sets can be easily derived using sequential arguments has also been, no doubt, a source of relief to the interested mathematics instructor.

Recently, Cakalli  \cite{Cakallisequentialdefinitionsofconnectedness} has defined $G$-sequential connectedness of a topological group and investigated some results in this generalized setting.

The purpose of this paper is to develop some further properties of $G$-sequential connectedness in metrizable topological groups, and present some interesting results.

\maketitle

\section{Preliminaries}

Before   giving  some results on $G$-sequential connectedness  we remark  some background as follows. Throughout this paper, $\textbf{N}$ denotes the set of all positive integers and  $X$  denotes a  Hausdorff topological group written additively satisfying  the first axiom of countability. We use boldface letters $\bf{x}$, $\bf{y}$, $\bf{z}$, ... for sequences $\textbf{x}=(x_{n})$, $\textbf{y}=(y_{n})$, $\textbf{z}=(z_{n})$, ... of terms of $X$. $s(X)$ and $c(X)$ respectively denote the set of all $X$-valued sequences and the set of all $X$-valued convergent sequences of points in $X$.

Following the idea given in a 1946 American Mathematical Monthly problem \cite{Buck}, a number of authors Posner \cite{Posner}, Iwinski \cite{Iwinski},
Srinivasan \cite{Srinivasan}, Antoni \cite{Antoni}, Antoni and Salat \cite{AntoniandSalat}, Spigel and Krupnik \cite{SpielandKrupnik} have studied $A$-continuity defined by a regular summability matrix $A$. Some authors \"{O}zt\"{u}rk \cite{Ozturk}, Sava\c{s} \cite{Savas}, Sava\c{s} and Das \cite{SavasandDas}, Borsik and Salat \cite{BorsikandSalat} have studied $A$-continuity for methods of almost convergence or for related methods. See also \cite{Boos} for an introduction to  summability matrices.

A sequence $(x_{k})$ of points in $X$ is called to be statistically convergent to an element $\ell$ of $X$ if for each neighborhood $U$ of $0$
\[
\lim_{n\rightarrow\infty}\frac{1}{n}|\{k\leq n: x_{k}-\ell \notin U\}|=0,
\] and this is denoted by $st-\lim_{n\rightarrow\infty}x_{n}=\ell$. Statistical limit is an additive function on the group of statistically convergent sequences of points in $X$ (See \cite{Fridy} for the real case and \cite{Cakallilacunarystatisticalconvergenceintopgroups}, \cite{Cakalli2}, \cite{CakalliAstudyonstatisticalconvergence} for the topological group setting and see \cite{MaioKocinacStatisticalconvergenceintopology}, and \cite{CakalliandKhan} for the most general case, i.e., topological space setting).

A sequence $(x_{k})$ of points in $X$ is called lacunary statistically convergent to an element $\ell$ of $X$ if
\[
\lim_{r\rightarrow\infty}\frac{1}{h_{r}}|\{k\in I_{r}: x_{k}-\ell \notin U\}|=0,
\]
for every neighborhood $U$ of 0 where $I_{r}=(k_{r-1},k_{r}]$ and $k_{0}=0$, $h_{r}:k_{r}-k_{r-1}\rightarrow
\infty$ as $r\rightarrow\infty$ and $\theta=(k_{r})$ is an increasing sequence of positive integers.
For a constant lacunary sequence $\theta=(k_{r})$, the lacunary statistically convergent sequences in a topological group form a subgroup of the group of all $X$-valued sequences and lacunary statistical limit is an additive function on this space (see \cite{Cakallilacunarystatisticalconvergenceintopgroups} for topological group setting and see \cite{FridyandOrhan1} and \cite{FridyandOrhan2} for the real case). Throughout this paper, we assume that  $\liminf_{r}\frac{k_{r}}{k_{r-1}}>1$.

By a method of sequential convergence, or briefly a method, we mean an additive function $G$ defined on a subgroup $c_{G}(X)$ of $s(X)$ into $X$
\cite{CakalliSequentialdefinitionsofcompactness}. A sequence \; \; $\textbf{x}=(x_{n})$ is said to be $G$-convergent to $\ell$ if $\textbf{x}\in c_{G}(X)$ and $G(\textbf{x})=\ell$. In particular, $\lim$ denotes the limit function $\lim \textbf{x}=\lim_{n}x_{n}$ on the group $c(X)$. A method $G$ is called {\em regular} if every convergent sequence $\textbf{x}=(x_{n})$ is $G$-convergent with $G(\textbf{x})=\lim \textbf{x}$.
Clearly if $f$ is $G$-sequentially continuous on $X$, then it is G-sequentially continuous on every subset $Z$ of $X$, but the converse is not necessarily true since in the latter case the sequences $\bf{x}$' s are restricted to $Z$. This was demonstrated by an example in \cite{ConnorGrosse} for a real function.

The notion of regularity introduced above coincides with the classical notion of regularity for summability matrices and with regularity in a topological group for limitation methods. See \cite{Boos} for an introduction to regular summability matrices, \cite{CakalliThorpe} for an introduction to regular limitation (summability) methods and see \cite{Zymund} for a general view of sequences of reals or complex.

 We recall the definition of $G$-sequential closure of a subset of $X$ from \cite{CakalliSequentialdefinitionsofcompactness} and  \cite{CakalliOnGcontinuity}. If  $A\subset X$ and $\ell \in X$, then  $\ell$ is in the $G$-sequential closure of $A$, which is called $G$-hull of $A$ in \cite{ConnorGrosse}, if there is a sequence $\textbf{x}=(x_{n})$ of points in $A$ such that $G(\textbf{x})=\ell$. Denote $G$-sequential closure of a set $A$ by $\overline{A}^{G}$  and say that a subset $A$ is {\em $G$-sequentially closed} if it contains all of the points in its $G$-closure, i.e.,  $\overline{A}^{G}\subset A$.

It is clear that $\overline{\phi}^{G}=\phi$ and $\overline{X}^{G}=X$ for a regular method $G$. If\; $G$ is a regular method, then $A\subset \overline{A}\subset \overline{A}^{G}$ and hence $A$ is $G$-sequentially closed if and only if $\overline{A}^{G}=A$. Note that even for regular methods, it is not always true that $\overline {(\overline{A}^{G})}^{G}=\overline{A}^{G}$.

Even for regular methods, the union of any two $G$-sequentially closed subsets of $X$ need not be a $G$-sequentially closed subset of $X$ as  seen  considering  Counterexample 1  in \cite{CakalliOnGcontinuity}.

\c{C}akall\i \;  has introduced the concept of $G$-sequential compactness and proved that $G$-sequentially continuous image of any $G$-sequentially compact subset of $X$ is also $G$-sequentially compact  \cite[Theorem 7] {CakalliSequentialdefinitionsofcompactness}. He has defined  $G$-sequential continuity and obtained further results in \cite{CakalliOnGcontinuity} (see also \cite{CakalliForwardcontinuity}, \cite{CakalliDeltaquasiCauchysequences}, \cite{Cakalli6}, \cite{DikandCanak}, \cite{CakCanDik}, \cite{CakalliNewkindsofcontinuities}, \cite{Cakallistatisticalwardcontinuity} for some other types of continuities which can not  be given by any sequential method and see \cite{CakallistatisticalquasiCauchysequences} for some kinds of continuities which coincide with uniform continuity when the domain of the function is connected). A function $f\colon X \rightarrow X$ is {\em $G$-sequentially continuous} at a point $u$ if, given a sequence $(x_{n})$ of points in $X$, $G(\textbf{x})=u$ implies that $G(f(\textbf{x}))=f(u)$.

Recently Mucuk and Sahan \cite{MucukSahan} have investigated further properties of $G$-sequential closed subsets of $X$.  They have modified the definition of open subset to $G$-sequential case in the sense that a subset $A$ of $X$  is {\em $G$-sequentially open} if its complement is $G$-sequentially closed, i.e., $\overline{X\setminus A}^{G}= X\setminus A$ and obtained that the union of any $G$-sequentially open subsets of $X$ is $G$-sequentially open. From the fact that for a regular sequential method $G$,   $G$-sequential closure of a subset of $X$ includes the set itself  it is easy to see that a subset $A$ is $G$-sequentially open if and only if  $\overline{X\setminus A}^{G}\subseteq X\setminus A$.  If a function $f$ is $G$-sequentially continuous on $X$, then inverse image $f^{-1}(K)$  of any $G$-sequentially closed subset $K$ of $X$ is $G$-sequentially closed \cite{CakalliSequentialdefinitionsofcompactness}. If a function $f$ is $G$-sequentially continuous on $X$, then inverse image $f^{-1}(U)$  of any $G$-sequentially open subset $U$ of $X$ is $G$-sequentially open \cite[Theorem 12]{MucukSahan}.  For a regular method $G$ function $f_a\colon X\rightarrow X,x\mapsto a+x$ is $G$-sequentially continuous, $G$-sequentially closed  and $G$-sequentially open \cite[Corollary 3]{MucukSahan} for any constant $a$ in $X$. If $A$ and $B$ are $G$-sequentially open, then so is the sum $A+B$ \cite[Theorem 15]{MucukSahan}. Mucuk and Sahan  have also obtained that a subset $A$ of $X$ is $G$-sequentially open if and only if  each $a\in A$ has a $G$-sequentially open neighborhood $U_a$ such that $U_a\subseteq A$ \cite[Theorem 4]{MucukSahan}.

Huang and Lin \cite{Huang-Lin}; and Fedeli and Donne \cite{Fedeli-Donne} have  investigated sequential connectedness.

Cakalli \cite{Cakallisequentialdefinitionsofconnectedness} has recently introduced the concept of $G$-sequential connectedness as follows.

\begin{Def} \label{Defseqconnected}   A non-empty subset $A$ of a topological group $X$ is called  {\em $G$-sequentially connected} if there are no non-empty,  disjoint $G$-sequentially closed subsets $U$ and $V$ of $X$ meeting $A$ such that  $A\subseteq U\bigcup V$.    Particularly  $X$ is called $G$-sequentially connected, if there are no non-empty,  disjoint $G$-sequentially closed subsets of $X$  whose union is $X$.
\end{Def}

The following Theorem is given  in   \cite[Corollary 1]{Cakallisequentialdefinitionsofconnectedness}.
 \begin{Thm}\label{Corclosurecon} If $G$ is a regular sequential method and $A$ is a $G$-sequentially connected subset of $X$, then so is $\overline{A}^G$. \end{Thm}

\begin{Thm} \label{TheUnionconn} \cite[Theorem 3]{Cakallisequentialdefinitionsofconnectedness} Let  $\{A_i\mid i\in I\}$ be a class of $G$-sequentially connected subsets of $X$ if $\bigcap_{i\in I}A_i$ is non-empty, then $\bigcup_{i\in I}A_i$ is $G$-sequentially connected. \end{Thm}

\begin{Thm} \label{Theoopensubgpclosed} \cite[Theorem 5]{Cakallisequentialdefinitionsofconnectedness}  Let  $G$ be a regular sequential method, and $H$ a subgroup of $X$. If $H$ is $G$-sequentially open, then it is $G$-sequentially closed.\end{Thm}

We recollect the following definition.

\begin{Def} \label{DefclosedsubsetinA} \cite[Definition 2]{Cakallisequentialdefinitionsofconnectedness} Let $A$ be a subset of $X$. A subset $F\subseteq A$ is called $G$-{\em  sequentially closed} in $A$ if  $F=U\cap A$ for some $G$-sequentially closed subset $U$ in $X$. We say that  a subset $V\subseteq A$ is  $G$-{\em sequentially open} in $A$ if $A\backslash V$ is $G$-sequentially closed in $A$.\qed \end{Def}

Here we remark that  a subset $B\subseteq A$ is $G$-sequentially open in $A$ if and only if $B=A\cap V$ for a $G$-sequentially open subset $V$ of $X$.

\section{Results}

First, we introduce a concept of $G$-sequentially connected component of a point $x$ in $X$ which extends the concept of an ordinary sequential connected component of a point $x$ in $X$.

  \begin{Def}
  The  largest $G$-sequentially connected subset containing a point $x$ in $X$ is called {\em $G$-sequentially connected component} of $x$ and denoted by  ${C_x}^G$ .\end{Def}
 We note that ${C_x}^G$ coincides with the ordinary sequential connected component of $x$ when $G=lim$. We write ${\pi_0X}^G$ for the set of $G$-sequentially connected components of all points in $X$ and similarly write ${\pi_0A}^G$ for the set of all $G$-sequentially connected components of all points in a subset $A$.

  \begin{Lem} \label{Inthesamecomp} Let  $x,y\in X$. If $x$ and $y$ are in a $G$-sequentially connected subset $A$ of $X$, then $x$ and $y$ are in the same $G$-sequentially component of $X$.\end{Lem}
\begin{Prf} Let $x$ and $y$ be contained in a $G$-sequentially connected subset $A$ of $X$. Then $x,y\in A\subseteq {C_x}^G$ and $x,y\in A\subseteq {C_y}^G$.  So ${C_x}^G\subseteq {C_y}^G$ and ${C_y}^G\subseteq {C_x}^G$. Therefore  $C_x= C_y$. \end{Prf}

   \begin{Lem}\label{Concomparti} The $G$-sequentially connected components of $X$ form a partition of $X$. \end{Lem}

  \begin{Prf} It is obvious that  $G$-sequentially connected components form a cover of $X$. We prove that for $x,y\in X$ if the components  ${C_x}^G$ and ${C_y}^G$ intersect, then ${C_x}^G={C_y}^G$.  Let $z\in {C_x}^G\cap {C_y}^G$. Then $z\in  {C_x}^G \subseteq {C_z}^G$ and  $z\in  {C_y}^G \subseteq {C_z}^G$ since  ${C_z}^G$ is the largest $G$-sequentially connected subset including $z$.  On the other hand  ${C_z}^G \subseteq {C_x}^G$ and ${C_z}^G \subseteq {C_y}^G$ since $x\in {C_x}^G\subseteq {C_z}^G$ and $ y\in {C_y}^G\subseteq  {C_z}^G$. Therefore ${C_x}^G={C_z}^G={C_y}^G$, and so  the proof is completed. \end{Prf}

 \begin{Thm}\label{Concompbijectif} Let $A,B\subseteq X$. If $A$ and $B$ are $G$-sequentially homeomorphic, then
 ${\pi_{0}A}^G$ and ${\pi_{0}B}^G$ have the same cardinality, i.e. there exists a bijection between them. \end{Thm}

\begin{Prf} Let $f\colon A\rightarrow B$ be a $G$-sequentially homeomorphism.  Define a map \[ \pi_0(f)\colon {\pi_{0}A}^G\rightarrow {\pi_{0}B}^G,~f({C_x}^G)={C_{f(x)}}^G\] induced by the function $f\colon A\rightarrow B$.  Since $f$ is $G$-sequentially continuous, the image of a $G$-sequentially connected subset is $G$-sequentially connected \cite[Theorem 1]{Cakallisequentialdefinitionsofconnectedness}.  Hence if $y\in {C_{x}}^G$ then by Lemma \ref{Inthesamecomp}, $f(x)$ and $f(y)$ are in the same $G$-sequentially connected component and so   ${C_{f(x)}}^G={C_{f(y)}}^G$. Therefore the  map $\pi_0(f)$ is well defined. Since $f^{-1}$ is $G$-sequentially continuous, if ${C_{f(x)}}^G=C_{f(y)}^G$, then $G$-sequentially components of $x$ and $y$ are same, i.e, ${C_{x}}^G={C_{y}}^G$ and therefore $\pi_0(f)$ is injective. Further since  $f({C_{x}}^G)={C_{y}}^G$  with  $y=f(x)$,  the map   $\pi_0(f)$ is onto. \end{Prf}
\begin{Thm} \label{Concomcon}  $G$-sequentially connected component of a point $x$ in $X$ is $G$-sequentially closed for any regular sequential method $G$.
\end{Thm}
\begin{Prf}  Since the $G$-sequentially connected  component $C_x^G$ is $G$-sequentially connected, by Theorem \ref{Corclosurecon}  the closure  $\overline{{C_x}}^G$ is $G$-sequentially connected.   Since $G$ is regular  ${C_x}^G\subseteq {\overline{C_x}}^G$ but the largest $G$-sequentially connected subset including $x$ is ${C_x}^G$. Therefore    $\overline{{C_x}}^G\subseteq {C_x}^G$ and so  ${C_x}^G$ is $G$-sequentially closed.
\end{Prf}
\begin{Thm} \label{LemConneccompoforigin} For a    regular sequential method $G$, the $G$-sequentially  connected component of the identity is a $G$-sequentially closed,  normal subgroup of $X$. \end{Thm}
\begin{Prf} Write   $K$ for the $G$-sequentially connected  component of the identity point $0$.
By Theorem \ref{Concomcon},  $K$ is $G$-sequentially closed.  To prove that $K$ is a subgroup, we prove that $K-K\subseteq K$, where $K-K$ is the set of all points $x-y$ for $x,y\in K$.  Since $K$ is a $G$-sequentially connected subset,  for each $x\in K$, the set    $x-K$ is $G$-sequentially connected as the image of  $G$-sequentially connected subset $K$ under a $G$-sequentially continuous function.  Then by Theorem  \ref{TheUnionconn}
\[K-K=\bigcup_{x\in K}(x-K)\] is $G$-sequentially connected as a union of $G$-sequentially connected subsets including  $0$. But the largest $G$-sequentially connected subset including $0$ is $K$.   Therefore  $K-K\subseteq K$, i.e, $K$ is a subgroup. Further for any $a\in X$ the function \[f_a\colon X\rightarrow X,x\mapsto a+x-a\] is $G$-sequentially continuous. So $f_a(K)=a+K-a$ is a $G$-sequentially connected subset  including  $0\in X$. Therefore $a+K-a\subseteq K $, which implies that $K$ is normal.
\end{Prf}

Let  $X$ be a topological group and  $U$ a symmetric neighbourhood of the identity $0$. We say that $X$ is {\em generated} by $U$, if each element of $X$ can be written as a sum of some elements in $U$.
\begin{Thm} Let  $G$ be a regular sequential method. If $X$ is generated by a $G$-sequentially connected  and symmetric neighbourhood $U$ of the identity, then $X$ is $G$-sequentially connected.\end{Thm}

\begin{Prf}  We prove that every point of $X$ is in the $G$-sequentially connected component of $0$. Let $x\in X$. Since $X$ is generated by $U$, the point $x$ can be written as $x=x_1+x_2\dots+x_n$ for some $x_i\in U$. So $x\in (U+U+\dots +U)$ which is $G$-sequentially connected by \cite[Lemma 7]{Cakallisequentialdefinitionsofconnectedness}  and includes the identity point $0$. Therefore by Lemma \ref{Inthesamecomp}, $x$ and  $0$ are in the same $G$-sequentially connected component. This proves that $X$ has only one $G$-sequentially connected component, i.e, $X$ is $G$-sequentially connected. \end{Prf}

In the following theorem we prove that a $G$-sequentially connected topological group is generated by any $G$-sequentially open neighbourhood of the identity point of $X$.
\begin{Thm} Let  $G$ be a regular sequential method. If $X$ is $G$-sequentially connected, then $X$ can be generated by any  $G$-sequentially open and symmetric  neighborhood of the identity.\end{Thm}

\begin{Prf} Let $U$ be a $G$-sequentially open and symmetric  neighbourhood of the identity point $0\in X$ and let $X_U$ be the subgroup of $X$ generated by $U$.  Hence each $x\in X_U$ can be written as a sum $x=x_1+\dots +x_n$ of the points in $U$ and so $x\in (U+U+\dots +U)$. Therefore $X_U$ can be written as a union of the sums $U+U+\dots +U$. Here by \cite[Theorem 15]{MucukSahan} each sum  $U+U+\dots +U$ is  $G$-sequentially open and includes $0$.  By
\cite[Theorem 3]{MucukSahan}   $X_U$ becomes $G$-sequentially open and therefore by Theorem \ref{Theoopensubgpclosed},  $X_U$ is also $G$-sequentially closed. Since $X$ is $G$-sequentially   connected it follows that  $X_U=X$. \end{Prf}

In the following definition of $G$-sequentially locally connectedness by a $G$-sequentially connected neighbourhood of a point  $x\in X$, we mean a  $G$-sequentially connected subset which contains a $G$-sequentially  open  subset including  $x$.
\begin{Def} \label{DefLocGseqconnected}
Let $G$ be a sequential method on $X$. We call $X$ as {\em $G$-sequentially locally connected}, if for any  $G$-sequentially  open neighbourhood $U$ of $x$, there is a $G$-sequentially  connected neighbourhood $V$  of $x$ such that $x\in V\subseteq U$.\end{Def}

\begin{Thm}\label{PropGLocseqcon}   $X$ is $G$-sequentially locally connected if and only if  $G$-sequentially  connected components  of any  $G$-sequentially   open subset are $G$-sequentially open . \end{Thm}
\begin{Prf} Let $X$ be  $G$-sequentially locally connected. Let  $A$ be a $G$-sequentially  open subset of $X$,  $C$ a $G$-sequentially  connected component of $A$ and $x\in C$.  Since $X$ is $G$-sequentially locally connected, there is a $G$-sequentially  connected neighbourhoof $U_x$ of $x$ such that $U_x\subseteq A$. But since the largest $G$-sequentially  connected subset of $A$ containing $x$ is $C$,  we have that $x\in U_x\subseteq C$.  Therefore by \cite[Theorem 4]{MucukSahan}, $C$ is $G$-sequentially open.

On the other hand if  $G$-sequentially connected components  of  any $G$-sequentially open subset is $G$-sequentially open, then $X$ becomes $G$-sequentially locally connected.\end{Prf}

A special case of Theorem \ref{PropGLocseqcon} is that if $X$ is $G$-sequentially locally connected, then each $G$-sequentially connected component of $X$ is $G$-sequentially open.

\begin{Thm}\label{PropGLocseqconn} Let $G$ be a regular method on $X$ and  $A,B\subseteq X$. Let $f\colon A\rightarrow B$ be an onto, $G$-sequentially continuous and $G$-sequentially open function. If $A$ is  $G$-sequentially locally connected, then so is $B$.  \end{Thm}
\begin{Prf} Suppose that  $f\colon A\rightarrow B$ is  an onto  function which is $G$-sequentially continuous and $G$-sequentially  open. Let $a\in A$ and $b\in B$ such that $f(a)=b$, and let $U$ be a $G$-sequentially open neighbourhood of $b$ in $B$.  Since $f$ is $G$-sequentially continuous, by \cite[Theorem 12]{MucukSahan}  $f^{-1}(U)$ is a $G$-sequentially open neighbourhood of $a$. Since $A$ is locally  $G$-sequentially connected, there is a $G$-sequentially connected neighbourhood of $a$ such that $V\subseteq f^{-1}(U)$. This implies that  $f(V)\subseteq U$.  Here since $f$ is $G$-sequentially open, $f(V)$ is $G$-sequentially open  and since $f$ is $G$-sequentially continuous, $f(V)$ is $G$-sequentially connected. Therefore $B$ is also $G$-sequentially locally connected. \end{Prf}

\section{Conclusion}

The present work contains further results on $G$-sequentially connectedness in first countable and  Hausdorff topological groups. So that one may expect it to be more useful tool in the field of topology in modeling various problems occurring in many areas of science, geographic information systems, population modeling and motion planning in robotics. It seems that an investigation of the present work taking `nets' instead of `sequences' could be done using the properties of 'nets' instead of using the properties of `sequences'. As the vector space operations, namely, vector addition and scalar multiplication,  are continuous in a cone normed space so cone normed spaces are special topological groups, we see that the results are also valid in cone normed spaces (see \cite{SonmezandCakalli} for the definition of a cone normed space). For further study, we also suggest to investigate the present work for the fuzzy case. However, due to the change in settings, the definitions and methods of proofs will not always be analogous to those of the present work (see \cite{CakalliandPratul} for the definitions in the fuzzy setting).

\end{document}